\documentclass{article}
\usepackage{amsthm}
\usepackage{amsfonts}
\usepackage{cite}
\usepackage{amsmath, amscd}
\newtheorem{theorem}{Theorem}[section]
\newtheorem{definition}[theorem]{Definition}
\newtheorem{proposition}[theorem]{Proposition}
\newtheorem{lemma}[theorem]{Lemma}
\newtheorem{example}[theorem]{Example}
\newtheorem{corollary}[theorem]{Corollary}

\begin{document}
\author{Jeremy Berquist}
\title{On Semirational Singularities}
\maketitle

\noindent
\linebreak
\textbf{Abstract.  }  We study semiresolutions of quasi-projective varieties with properties $G_1$, $S_2$, and seminormality.  Equivalently, these are varieties $X$ with Serre's $S_2$ property, such that there exists an open subvariety $U$, with complement of codimension at least two, such that the only singularities of $U$ are (analytically) double normal crossings.  Such varieties have been called ``demi-normal" by Koll\'ar \cite{Kol13}.  First, we discuss why these properties are ideal for the study of nonnormal varieties that appear in the birational classification of varieties.  We define semiresolutions and provide examples to illustrate the procedure of gluing along the conductor as the fundamental tool in obtaining a semiresolution of $X$ from a resolution of its normalization.  As an application of these methods, we discuss semirational singularities.  Our main results are a semismooth Grauert-Riemenschneider vanishing theorem (Theorem 4.2) and a proof that the definition of semirational singularities does not depend on the semiresolution chosen (Theorem 4.3) .  The smooth version of (4.2) first appeared in \cite{GR70}.   We also explain why semirational singularities are, in particular, Cohen-Macaulay and DuBois.  Finally, we discuss the role semiresolutions play in forming a nonnormal interpretation of the results of de Fernex and Hacon in \cite{dFH09}, where the $\mathbb{Q}$-Cartier hypothesis was found to be extraneous in the birational classification of normal varieties.

\noindent
\linebreak
\textbf{Acknowledgement.  }  The author would like to thank S\'andor Kov\'acs for his generous help with the more complicated pieces of this article, including the proof of semismooth G-R vanishing.  The results here appeared in the author's dissertation ``Singularities on Nonnormal Varieties," which was written under the advisement of Professor Kov\'acs.
\newpage
\tableofcontents
\addcontentsline{}{}{}
\newpage
\begin{section}{Introduction}
We have often seen in the study of varieties with nice properties, that it becomes necessary to study varieties without those properties.  A family of irreducible varieties may degenerate to a reducible variety, for instance.  Or, classifying varieties of general type may involve studying singularities, even when starting with a nonsingular variety.  In many cases, the object of study is a smooth variety paired with a simple normal crossing divisor.  This divisor, viewed as a scheme, has special singularities.  In such a situation, it might be necessary to resolve the singularities of the divisor in the best possible way.

Resolution of singularities is a powerful technique in studying the singularities of a given variety $X$.  In his monumental 1964 paper, Hironaka proved that resolution of singularities is possible for all varieties over a field of characteristic zero.  In fact, all of the methods of semiresolution are based on the existence of resolutions.  For normal varieties, we are fortunate enough to be able to keep the codimension one part of the variety intact.  In particular, we can compare divisors on $X$ with divisors on the resolution.  However, for a nonnormal variety (with the properties we discuss below), a resolution is too robust in this sense.  There is no such morphism that keeps the codimension one part intact when $X$ is nonnormal.  See (3.2) below.

The typical model for the varieties we consider here are simple normal crossing divisors (in an ambient smooth variety).  The components of such a scheme are smooth and intersect everywhere transversally, and in codimension one these intersections are double normal crossings.  For more general purposes, it is the behavior in codimension one that we would like to retain, and not the behavior in larger codimensions.  The properties we consider on $X$ are that it has Serre's condition $S_2$, and that there is an open subvariety $U$ whose complement has codimension at least two in $X$, such that the only singularities in $U$ are double normal crossings.  That is, at a closed singular point $x$, the completion of the local ring, $\hat{\mathcal{O}}_{X,x}$, is isomorphic to $k[[x_1, x_2, \ldots, x_n]]/(x_1x_2).$  Equivalent conditions on $X$ are that it has the properties $S_2$, $G_1$ (Gorenstein in codimension one), and seminormality.  (Note that a normal variety is one with properties $S_2$ and $R_1$, which is slightly stronger than $G_1$, and is seminormal.)  

With such varieties as these, we must also choose the appropriate analog of a resolution of singularities $f: Y \rightarrow X$.  We cannot ask for $Y$ to be nonsingular.  However, if we require only that $Y$ is semismooth, then we can obtain a proper, birational morphism that is an isomorphism over the open set $U$ of $X$ described above, and such that the other singularities of $Y$ are relatively mild.  We call $Y$ semismooth if its closed points are either smooth points, double normal crossing points, or pinch points.  At a pinch point $x$, the completion of the local ring is isomorphic to $k[[x_1, x_2, \ldots, x_n]]/(x_1^2-x_2^2x_3).$  That this type of singularity should appear along with double normal crossings is due in part to the way that nonnormal varieties are obtained from their normalizations.  We will call a morphism $f: Y \rightarrow X$ a semiresolution if it is proper (we can even ask that it be projective), birational, an isomorphism over the open subvariety $U \subseteq X$ described above, such that $Y$ is semismooth and no component of the singular locus of $Y$ is exceptional.

Varieties with these properties share several similarities with their normal counterparts.  A theory of divisors can be developed with just the conditions $G_1$ and $S_2$.  See \cite{Hart94}.  In addition, there is a semismooth version of the Grauert-Riemenschneider vanishing theorem, which originally states that for a resolution of singularities $f: Y \rightarrow X$ (more generally, a proper birational morphism from a smooth variety to $X$), $R^if_*\omega_Y = 0$ for $i >0.$  That this condition holds for a semiresolution of a variety with the above properties is shown in Theorem 4.2.  One may also define semirational singularities in a manner that parallels the definition of rational singularities.  These are defined as normal varieties $X$ such that there exists a proper, birational morphism $g: Y \rightarrow X$ from a nonsingular variety $Y$ to $X$, and where $R^ig_*\mathcal{O}_Y = 0$ for $i > 0$.  It should be checked that this definition does not depend on the choice of $g$; then one can choose a resolution of singularities as the defining morphism.  For nonnormal varieties, it is also true that if this condition holds for one semiresolution, then it holds for all semiresolutions.  We prove this last fact in Theorem 4.3.  (Note that these results are stated for varieties over a field of characteristic zero.  We assume this throughout the present article.)

The paper is organized as follows.  First, we review the definitions of $S_2$, $G_1$, and seminormality and present the basic definitions needed to understand what semiresolutions are.  In the next section, we present several examples.  These are meant to illustrate the basic tool of gluing along the conductor in order to construct a semiresolution of $X$ from a resolution of its normalization.  Finally, as an application of the technique of passing to the normalization, we define semirational singularities and prove a semi-smooth version of G-R vanishing.  We observe that semirational varieties are Cohen-Macaulay and DuBois, and give an example of a nonnormal variety without semirational singularities.  

All varieties are assumed to be reduced and over an algebraically closed field of characteristic zero.  In particular, a variety may be reducible, and resolution of singularities is true for such varieties.  Also, since we are concerned mostly with local properties, a variety is assumed to be quasi-projective unless otherwise mentioned.
\end{section}

\begin{section}{Preliminaries}
We first recall the definitions of $S_2$, $G_1$, and SN.

\begin{definition}  A coherent sheaf $\mathcal{F}$ on a variety $X$ is $S_n$ provided that for all points $x$, we have $$\textnormal{depth } \mathcal{F}_x \geq \textnormal{min}(\textnormal{dim } \mathcal{F}_x, n).$$
\end{definition}

We remark that some authors use a similar definition involving the dimension of the local ring $\mathcal{O}_{X,x}$ in place of the dimension of the module $\mathcal{F}_x$.  The difference between the two definitions is in whether one wants to think of a module as a module over a ring or over the quotient via the annihilator of the module.

A coherent sheaf is \textit{Cohen-Macaulay} if it is $S_n$ for all $n$.  Thus, considering that the depth of a module is always bounded above by its dimension, a Cohen-Macaulay sheaf is one for which the depth and dimension are the same at each localization.

Of course, a variety will be called $S_n$ (or Cohen-Macaulay) if its structure sheaf $\mathcal{O}_X$ has this property.

\begin{definition}  A variety is said to have condition $G_1$ if $\mathcal{O}_{X,x}$ is a Gorenstein ring whenever $x$ is a point of codimension 0 or 1.  In other words, at such points, the canonical module is trivial and the variety is Cohen-Macaulay.
\end{definition} 
\noindent
The canonical module does not have to be defined here; we just note that the dualizing sheaf exists for all of our varieties, and that the localizations of the dualizing sheaves are the canonical modules for the corresponding local rings.  Thus, $G_1$ for us simply means that the dualizing sheaf is invertible in codimensions 0 and 1.

We have the following implications for local rings:  $$\textnormal{regular} \implies \textnormal{Gorenstein} \implies \textnormal{Cohen-Macaulay}.$$

\noindent
In fact, any complete intersection is Gorenstein.  In particular, hypersurfaces are Gorenstein, hence $G_1$ and $S_2$.

\begin{definition}  An extension of rings $A \hookrightarrow B$ is a quasi-isomorphism if it is finite, a bijection on prime spectra (hence a homeomorphism on prime spectra), and each residue field extension $k(p) \hookrightarrow k(q)$ is an isomorphism.
\end{definition}
\noindent
The term ``subintegral" is also used in place of ``quasi-isomorphism."

\begin{definition}  Given a finite extension of rings $A \hookrightarrow B$, the seminormalization of $A$ in $B$ is the unique largest subring of $B$ that is quasi-isomorphic to $A$.  We say that $A$ is seminormal in $B$ if it equals its seminormalization in $B$.  We say that a reduced ring $A$ is seminormal if it equals its seminormalization in its integral closure.
\end{definition}

We note that the normalization is finite for reduced finitely-generated algebras over a field.  It is not hard to show that $A$ is seminormal in $B$ if $b^2, b^3 \in A$ imply $b \in A$ for any $b \in B.$

We say that a (reduced) variety is seminormal if each of its local rings is seminormal.

\begin{proposition}[\cite{Kol13}, 5.1] A variety $X$ is $S_2$, $G_1$, and SN if and only if it is $S_2$ and there exists an open subvariety $U$ such that:  (i)  the complement of $U$ in $X$ has codimension at least two; (ii)  for any closed singular point $x \in U$, $\hat{\mathcal{O}}_{X,x} \cong k[[x_1, x_2, \ldots, x_n]]/(x_1x_2).$  
\end{proposition}

Why should we care about these varieties?  Varieties such as these arise in several situations.  Obviously, a normal variety has these conditions.  As another example, a strong resolution of singularities of a variety $X$ produces a nonsingular variety such that the preimage of the singular locus of $X$ is a divisor with simple normal crossings.  This divisor, viewed as a variety in its own right, is $S_2$ as a hypersurface and clearly has only double normal crossings in codimension one.

The two conditions $S_2$ and $G_1$ together constitute slightly weaker conditions than normality and for which a theory of generalized divisors can be worked out. See \cite{Hart94}.  The irreducible codimension one subschemes correspond to reflexive sheaves that are invertible at the generic points of $X$.  Thus one can work on the level of sheaves, which are in some ways more flexible than divisors.

The conditions $S_2$ and SN together imply that the nonnormal locus of $X$ is a reduced divisor.  Given a reduced ring $A$ with finite normalization $B$, the conductor $(A:B) := \{a \in A: aB \subseteq A\}$ is the largest ideal of $A$ that is simultaneously and ideal in $B$.  With $S_2$ and SN, the conductor is reduced and has as associated points only height one primes.  When glued properly, these ideals become the divisor in $X$ that defines the nonnormal locus.

For the study of singularities on nonnormal varieties, we thus have an intricate relationship between these three properties.  They are, in a sense, the weakest possible conditions for such a study.

We now come to the definition of a semiresolution.  Recall that a pinch point is a point whose local ring is analytically isomorphic to $k[[x_1, x_2, \ldots, x_n]]/(x_1^2-x_2^2x_3).$  The pinch point is significant in that blowing up the origin produces another pinch point, and thus cannot be simplified without blowing up the entire double locus.  A pinch point on a surface is a quotient of a double normal crossing point.  See [7, 10.4] for more details on the relationship between double normal crossings and pinch points.

\begin{definition}
A variety is semismooth if every closed point is either smooth, a double normal crossing point, or a pinch point.
\end{definition}

\begin{definition}  Let $X$ be a variety as above.  Then a morphism $f: Y \rightarrow X$ is a semiresolution if the following conditions are satisfied:  (i)  $f$ is projective, (ii) $Y$ is semismooth, (iii) no component of the conductor $C_Y$ is $f$-exceptional, and (iv) $f|_{f^{-1}(U)}: f^{-1}(U) \rightarrow U$ is an isomorphism, where $U$ is an open set whose only closed singular points are double normal crossings.
\end{definition}

If $U \subseteq X$ is an open subvariety whose closed singular points are (analytically) double normal crossing points or pinch points, then Koll\'ar shows that there is a projective morphism $f: Y \rightarrow X$ such that $Y$ is semismooth, $f^{-1}(U) \rightarrow U$ is an isomorphism, and the singular (double) locus of $Y$ maps birationally onto the singular locus of $U$.  Thus semiresolutions exist.  We should stress that the characteristic of the base field is zero; Koll\'ar's construction uses a resolution of singularities of the normalization of $X$.

There is also an analog of log resolution of singularities for nonnormal varieties.  With $S_2$ and $G_1$, we can speak of Weil divisors on $X$ as formal sums of codimension one subschemes, and these correspond to rank one reflexive sheaves.  Although the notion of pullback of Weil divisors is somewhat delicate (there are issues with the group laws and with compatibility with intersection theory), pulling back reflexive sheaves is not difficult.  It is possible to pull back sheaves in a well-defined way, so that the resulting theory makes sense.  See \cite{dFH09} for the normal version.  The author found the nonnormal version in his dissertation.

In the present case, we want the exceptional divisor and the conductor to be distinct and to have simple normal crossings.  But, for a semismooth variety $Y$, there is not a local system of parameters, so we need to be careful when we talk about simple normal crossing divisors.  As we see in the next section, the normalization $\overline{Y}$ is smooth.  The correct interpretation is for the exceptional set and the conductor to lift to a simple normal crossing divisor on $\overline{Y}$.  Equivalently, before lifting, we have a scheme with smooth components that intersect transversally.  

Both double normal crossing points and pinch points are hypersurface singularities.  Locally analytically (even locally), $X$ behaves like a subvariety of some affine space $\mathbb{A}^n$.  One can form a simple normal crossing divisor $D$ in $\mathbb{A}^n$ that does not contain $X$, and restrict to obtain a divisor on $X$.  In case $X$ is locally analytically a double normal crossing point or pinch point, the conductor is a smooth subscheme of $\mathbb{A}^n$, and it can be asked whether it has transverse intersections with the given divisor.  Thus, by restricting, we obtain the general model for a simple normal crossing divisor on a semismooth variety.

The local models for such a divisor $D'$ are therefore as follows:  
\begin{itemize}
\item{Near a smooth point, $X$ is given by $x_1 = 0$ and $D' = (\Pi_{i \in I}x_i = 0)$ for $I \subset \{2, \ldots, n+1\}$}
\item{Near a double normal crossing point, $X$ is given by $(x_1x_2 = 0) \subset \mathbb{A}^n_{x_1, \ldots, x_n}$, and the local model for $D'$ is $D' = (\Pi_{i \in I} x_i = 0)$, for some set $I \subset \{ 3, \ldots, n \}$. }
 \item{Near a pinch point, given locally analytically by $(x_1^2 - x_2^2x_3 = 0)$, the local model is $D' = (\Pi_{i \in I} x_i = 0) + D_2$, for some $I \subset \{ 4, \ldots, n \}$ and where either $D_2 = 0$ or $D_2 = (x_1 = x_3 = 0)$.  }
\end{itemize}

\noindent
In both cases, the lifting of $D$ to the normalization has simple normal crossings, and intersects the double locus transversally.  Likewise, if we start with a simple normal crossing divisor that contains the conductor as a component, we then glue along the conductor to obtain $D'$.  We will discuss the gluing construction in the next section.

In order to describe the nonnormal equivalent of a log resolution of a pair $(X,D)$, we must take into account the singularities of $D$ as well. The following proposition should be taken as the definition of a log resolution of singularities for nonnormal varieties.  In particular, one must prepare a divisor $D$ before attempting to apply a semi log resolution.

\begin{proposition}[\cite{Kol13}. 10.56] Let $X$ be a variety (seminormal with $G_1$ and $S_2$).  Suppose that $U$ is an open subset with only double normal crossing singularities.  Let $D$ be a Weil divisor such that $D|_U$ is smooth and disjoint from $\textnormal{Sing}(U)$.  Then there is a projective birational morphism $f: Y \rightarrow X$ such that $Y$ is semismooth, $f$ is an isomorphism over $U$, $\textnormal{Sing}(Y)$ maps birationally onto the closure of $\textnormal{Sing}(U)$, and the local models for $(Y, D' := f^{-1}_*D + \textnormal{Ex}(f))$ are as above.
\end{proposition}

\noindent
In fact, as ([7],10.56) shows, one obtains via $f$ a log resolution of $(\overline{X}, \overline{B}+ \overline{D})$, where $\overline{B}$ is the closure of the conductor of $f^{-1}(U)$ over $U$, and $\overline{D}$ is the preimage of $D$ in the normalization $\overline{X}$ of $X$.

We stress that the conditions on $D$ imply that the double locus has transverse intersections with $f^{-1}_*D + \textnormal{Ex}(f)$.  We include the double locus in our definition of normal crossings because it helps us keep track of the gluings we use to construct $Y$ and the total transform of $D$ from the normalization.  

We have now presented all the necessary definitions and are ready to start working through some examples.  We do so in the following section.  
\end{section}

\begin{section}{Examples of Semiresolutions}

We mentioned before that resolution of singularities is too robust a technique for studying the singularities of a nonnormal variety $X$.  In a bit more detail, what happens is that if $f: Y \rightarrow X$ is a proper birational morphism between nonnormal varieties, then $f$ is not automatically an isomorphism over the codimension one points of $X$.  This is the condition we need in order to compare divisors on $X$ with those on $Y$.  For varieties with $S_2$, this is equivalent to a condition on the structure sheaves.

\begin{example}  Let $R = k[x,y,z,w]/((x,y) \cap (z,w))$ be the coordinate ring of two planes in $\mathbb{A}^4$ meeting at a point.  Then the normalization $f:  Y \rightarrow X := \textnormal{Spec}(R)$ is an isomorphism over codimension one points, but $f_*\mathcal{O}_Y \neq \mathcal{O}_X$.
\end{example}

\noindent
The problem here is that $X$ is not $S_2$.  (In fact, at the origin, the element $x-z$ forms a maximal regular sequence.)  If we require $S_2$, then we get something better.  Note that the converse statement of the following lemma is true without any $S_2$ requirement.

\begin{lemma}  Suppose $f: Y \rightarrow X$ is a proper birational morphism between (reduced) varieties, where $X$ is $S_2$.  Then $f$ is an isomorphism over codimension one points of $X$ if and only if $f_*\mathcal{O}_Y = \mathcal{O}_X$.
\proof  Suppose that the codimension one points of $X$ are in the isomorphism locus, and consider the short exact sequence $$0 \rightarrow \mathcal{O}_X \rightarrow f_*\mathcal{O}_Y \rightarrow Q \rightarrow 0.$$  Injectivity follows from the fact that $f$ is birational and $\mathcal{O}_X$ is $S_1$.  By hypothesis, $Q$ has support in codimension at least two.  If $Q$ were nonzero, we could localize the above exact sequence at an associated point $p$ of $Q$.  Notice that since $f$ is proper, $f_*\mathcal{O}_Y$ is coherent, and therefore $(f_*\mathcal{O}_Y)_p$ is a finite $\mathcal{O}_{X,p}$-module.  Since $Y$ is reduced, $(f_*\mathcal{O}_Y)_p$ is $S_1$ as a ring.  It follows by finiteness that it is $S_1$ as a module over $\mathcal{O}_{X,p}$.  Consider the long exact sequence in cohomology $$0 \rightarrow \textnormal{H}^0_p(X, \mathcal{O}_{X,p}) \rightarrow \textnormal{H}^0_p(X, (f_*\mathcal{O}_Y)_p) \rightarrow \textnormal{H}^0_p(X, Q_p) \rightarrow \textnormal{H}^1_p(X, \mathcal{O}_{X,p}) \rightarrow \cdots.$$  Since the second and fourth (nonzero) terms are 0, by the above arguments, so is the third.  This is a contradiction, since at an associated point a module has nonvanishing local cohomology in dimension 0.

Suppose conversely that $f$ is a proper birational morphism with the property that $f_*\mathcal{O}_Y = \mathcal{O}_X$.  To reach a contradiction, assume that $p$ is a codimension one point not in the isomorphism locus of $X$.  Since $f$ is birational, there are only finitely many codimension one points with this property.  Thus there exists a neighborhood $p \in V$ such that $f^{-1}(V) \rightarrow V$ is an isomorphism at all codimension one points of $V$ except for $p$.  By shrinking $V$ further to avoid generic points of other components of the nonisomorphism locus (necessarily of codimension at least two), we may assume that $\{ p \}^-$ is exactly the set of points over which $f|_{f^{-1}(V)}$ is not an isomorphism.  

The condition on structure sheaves implies that $f$ has connected fibers.  Thus there must be exactly one codimension one point of $Y$ lying over $p$.  Call this point $q$.  Looking at the preimage of $\{ p\}^-$, there may be other components besides $\{ q \}^-$, and the images of these components may be removed downstairs in $X$.  In summary, over a suitable open set in $X$, we may assume that $\{ p\}^-$ is exactly the set of points over which $f$ is not an isomorphism, and that its preimage is $\{q\}^-$.

Give $\{q\}^-$ its reduced induced structure.  Since the dimensions of the integral varieties $\{p\}^-$ and $\{q\}^-$ are the same, we may apply (\cite{Hart87} II.Ex.3.7) to show that there is an open neighborhood of $p$ in $\{p\}^-$ over which $f$ is finite.  This open set is the restriction of some open set $U$ in $X$.  By the above remarks, $f: f^{-1}(U) \rightarrow U$ is quasi-finite everywhere.  Since the fibers of $f$ are connected, this restricted morphism is in fact injective.  Being proper and birational, it is a homeomorphism, and the condition on structure sheaves implies that it is an isomorphism.  This contradicts our choice of $p$.  Thus $f$ is in fact an isomorphism over all codimension one points of $X$.\qed
\end{lemma}

The normal codimension one points of a reduced variety are in the nonsingular locus.  The set of nonnormal codimension one points must therefore be finite, as they determine irreducible branches of the singular locus (a reduced variety is $R_0$).  For seminormal varieties with $G_1$ and $S_2$, the set of nonnormal codimension one points form a reduced divisor (in other words, the conductor has no embedded points).  In the most important cases, $X$ can be constructed from its normalization by gluing along the conductor.  According to (\cite{Art70} 3.1), the universal pushout can be used to glue components of the conductor without affecting the rest of the normalization of $X$.  Specifically, given a closed subscheme $B \hookrightarrow Y$ and a finite morphism $q: B \rightarrow B/\tau$, the universal pushout $$\begin{CD}
B           @>>>   Y  \\        
@VqVV               @VfVV  \\         
B/\tau      @>>>   Y' 
\end{CD}$$
has the property that $f: Y \rightarrow Y'$ is proper, agrees with $B \rightarrow B/\tau$ on $B$ and is an isomorphism elsewhere.  

Recall that the two singularities that we are willing to accept are double normal crossings ${\mathcal{O}_{X,x}}^* \cong k[[x_1, \ldots, x_n]]/(x_1x_2)$ and pinch points ${\mathcal{O}_{X,x}}^* \cong k[[x_1, \ldots, x_n]]/(x_1^2 - x_2^2x_3)$.  Each of these singularities can be constructed by gluing along the conductor.  We include the following two examples because they describe the normalizations of double normal crossings and pinch points, as well as showing how the gluing is done.

\begin{example}\end{example}  We show how the double normal crossing point is obtained via a gluing operation.  The normalization $\overline{X}$ is smooth, and it has two components, given analytically by its coordinate ring $$k[[x_2, x_3, \ldots, x_n]] \times k[[x_1, x_3, \ldots, x_n]].$$  The conductor $\overline{C}$ is given by the disjoint union of two hyperplanes $(x_2 = 0)$ and $(x_1 = 0)$ in the normalization, and by the ideal $(x_1,x_2)$ in $k[[x_1, \ldots, x_n]]/(x_1x_2)$.  Thus the conductor $C$ is a smooth divisor in $X$, and its lifting to $\overline{X}$ is a smooth double covering of $C$.  

If we identify the components of $\overline{C}$ along their common corresponding subvarieties, the universal pushout is the ring-theoretic pullback of the diagram 
$$\begin{CD}
R                      @>>>            k[[x_2, \ldots, x_n]] \times k[[x_1, x_3, \ldots, x_n]] \\
@VVV                                   @VVV                        \\                          
k[[x_3, \ldots, x_n]]  @>>>            k[[x_3, \ldots, x_n]] \times k[[x_3, \ldots, x_n]]      
\end{CD}$$  
In other words, we get all pairs $(p(x_2, \ldots, x_n),q(x_1, \ldots, x_n))$ such that $$p(0, x_3, \ldots, x_n) = q(0, x_3, \ldots, x_n).$$  These form a subring of the direct product that is isomorphic to the coordinate ring of the double normal crossing point.

\begin{example}\end{example}
Here we show how the pinch point is obtained from its normalization.  The pinch point's coordinate ring $k[[x_1, \ldots, x_n]]/(x_1^2 -x_2^2x_3)$ is isomorphic to the ring $k[[st,t,s^2, x_4, \ldots x_n]]$.  In these coordinates, the normalization $\overline{X}$ has coordinate ring $k[[s,t,x_4, \ldots, x_n]]$.  The conductor $C$ is given by the ideal $(st,t)$ and by the ideal $(t)$ in the normalization.  Thus $C$ and $\overline{C}$ are both smooth divisors, and $\overline{C}$ is a ramified double cover of $C$.  

There is a natural $\mathbb{Z}_2$-action on $\mathcal{O}_{\overline{C}} = k[[s,x_4, \ldots, x_n]]$, sending $s \mapsto -s$ and $x_i \mapsto x_i$ for all $i$.  If we glue along the quotient, the pullback diagram is 
$$\begin{CD}
R                         @>>> k[[s^2,x_4, \ldots, x_n]] \\
@VVV                           @VVV \\
k[[s,t,x_4, \ldots, x_n]] @>>> k[[s,x_4, \ldots, x_n]]
\end{CD}$$
The ring $R$ is a subring of $k[[s,t,x_4, \ldots, x_n]]$ consisting of power series of the form $f_0(s^2,x_4, \ldots, x_n) + tf_1(s,t,x_4, \ldots, x_n)$.  Thus we get back the coordinate ring of the pinch point.
\\ \\
Pinch points and double normal crossing points are the only singularities we allow on a semiresolution, and the morphism to $X$ is constructed as a universal pushout in all cases, so these examples are of fundamental importance.

The pinch point and double normal crossing point should probably occur together in any reasonable class of singularities containing either one.  Note that $(x_1^2-x_2^2x_3 = 0)$ has double normal crossings away from the origin, and a more complicated singularity at the origin.  Similarly, the pinch point can be obtained as a quotient of $(x_1x_2=0)$.  In two dimensions, for example, we can let $\mathbb{Z}_4$ act on $(xy=0) \subset \mathbb{A}^3_{x,y,z}$ by sending $x \mapsto -y, y \mapsto x, z \mapsto \sigma z$, where $\sigma$ is a primitive fourth root of unity.

\noindent
\linebreak
\textbf{Remark.  }  The above examples show that when $X$ is semismooth, its normalization $\overline{X}$ is smooth.  Moreover, its singular locus is a smooth divisor $D$, the lift to $\overline{X}$ of which is a smooth double cover of $D$, ramified along the pinch locus.  These facts are obvious for smooth closed points; however, a variety is smooth if it is smooth at its closed points, since localization preserves regularity, by (\cite{Hart87}, II.8.14A).  We also need to know:  (i) that normalization and completion commute, which does not happen for all rings, but for excellent rings like reduced finitely-generated algebras over a field, this is true; and (ii) that a local ring is regular if and only if its completion is a regular local ring.

\begin{example}  Let $k[x,y,z,w]/(x(y^2-z^2w))$ be the coordinate ring of an affine variety $X$.  Then the blowup of the ideal $(x,y,z)$ is a semiresolution.
\proof  This variety is a hypersurface singularity, hence it is Gorenstein and $S_2$.  Its components are seminormal, and the ideal sheaf of their intersection $(x,y^2-z^2w)$ is reduced.  Thus $X$ is seminormal by (\cite{LV81} 2.18).  There are pinch points in codimension one that we would like to keep intact when performing a semiresolution.  If we blow up the ideal $(x,y,z)$, then properties (i) and (iv) of (2.6) are automatic.  There are three charts on the blowup.  On the chart where $x=x',y=x'y',z=x'z',w=w'$, the blowup is a pinch point $(y'^2-z'^2w' = 0)$.  The double locus here is $(y'=z'=0)$, which is not exceptional.  On the chart where $x=x'y',y=y',z=y'z',w=w'$ the blowup is smooth.  Finally, if $x=x'z',y=y'z',z=z',w=w'$, then the blowup is $(x'(w'-y'^2)=0)$, which is a double normal crossing singularity.  The conductor is given by $(x'= w'-y'^2=0)$ and is not exceptional.  We conclude that the blowup is semismooth and no singular component is exceptional.\qed
\end{example} 

We will see that the property on the conductor/double locus of a semiresolution allows us to pass easily to the normalization.  Actually, the example above is somewhat special in that a semiresolution can be guessed at without knowing anything about gluing.  Semiresolutions are usually constructed by resolving the singularities of the normalization first, and then by gluing along the birational transform of the double locus (in a specified open subset $U \subset X$).  This property need not hold automatically, as the following example shows. 

\begin{example}  A morphism with properties (i), (ii), and (iv) of (2.7) and without property (iii).  \end{example}

Consider the hypersurface $X \subset \mathbb{A}^3$ defined by $(x^5-y^2z=0)$.  If we blow up the origin, then in the chart with coordinates $x=x', y=x'y', z= x'z'$, we obtain a pinch point $(x'^2-y'^2z'=0)$.  The exceptional divisor is given by $(x'=0)$ and therefore contains the entire double locus $(x'=y'=0)$.  In the chart with coordinates $x=x'y', y=y', z=y'z',$ we get a smooth variety.  In the remaining chart, we obtain the surface defined by the equation $x'^5z'^2 - y'^2 = 0.$  Away from the line $(x'=y'=0)$, there are only double normal crossings.

Now blow up $(x'=y'=0)$ in the third chart.  Since the double locus in the first chart is unaffected, there will still be a component of the conductor that is exceptional (for the composition of blowups).  We obtain a smooth chart and a chart with equation $x_2^3z_2^2-y_2^2 = 0.$  Again, away from the line $(x_2=y_2=0)$, we have only double normal crossings.  If we blow up this line, we obtain a smooth chart and a pinch point.  \qed

\begin{example}  Consider the surface defined by $k[x,y,z]/(xy)$, and let $\mathbb{Z}_2$ act on it by $x \mapsto -x, y \mapsto -y, z \mapsto -z$.  The resulting quotient is a nonnormal variety $X$ with two normal components.  A semiresolution is obtained by resolving the components and gluing along the birational transform of the line of intersection.
\end{example}  
By definition, the hyperquotient singularity $X$ is obtained as a residue class of a ring of invariants.  If $\mathbb{Z}_2$ acts on $k[x,y,z]$ as stated, the ring of invariants is $k[x^2,xy,y^2,xz,yz,z^2]$.  We get the coordinate ring of $X$ by annihilating the intersection of the ideal $(xy)$ with this ring.  The intersection ideal has generators $xy,x^2y^2,xyz^2,x^2yz,$ and $xy^2z$.  Therefore, in terms of generators and relations, the coordinate ring $\mathcal{O}_X$ is given by $$k[u_0,u_1,u_2,u_3,u_4,u_5]/(I + J),$$ where $$I = (u_0u_2-u_1^2, u_2u_5-u_4^2, u_0u_5-u_3^2,u_1u_3-u_0u_4, u_3u_4-u_1u_5)$$ and $$J = (u_1,u_0u_2,u_3u_4,u_0u_4,u_2u_3).$$  If we simplify, then we obtain $k[u_0,u_2,u_3,u_4,u_5]$ modulo the ideal $$(u_2u_5-u_4^2,u_0u_5-u_3^2,u_0u_2,u_3u_4,u_0u_4, u_2u_3).$$
Note that the spectrum of this ring has two components, given by $(u_0=u_3=0)$ and $(u_2=u_4=0)$.  Each of these defines a quadric cone, and the cones are identified along the line $\mathbb{A}^1 \cong \textnormal{Spec } k[u_5]$.  We can obtain the same ring as a pullback, where the maps $c$ and $d$ are the quotient maps:
$$\begin{CD}
\mathcal{O}_X                @>>>   k[u_0,u_3,u_5]/(u_0u_5-u_3^2)\\
@VVV                                @VcVV \\
k[u_2,u_4,u_5](u_2u_5-u_4^2) @>d>>  k[u_5]
\end{CD}$$  
If we compute the $\mathbb{Z}_2$-quotient of each component of $k[x,y,z]/(xy)$ and then glue, we get $\mathcal{O}_X$.  In fact, the quadric cone is a $\mathbb{Z}_2$-quotient of affine two-space $\mathbb{A}^2$.

The components of $X$ are both normal, and each of their singularities is resolved by blowing up the origin.  In each of these two blowups, the birational transform of the intersection line $\mathbb{A}^1$ appears in only one chart.  Gluing over the corresponding ring $k[u_5']$, we get the pullback diagram 
$$\begin{CD}
\mathcal{O}_Y @>>>    k[u_3',u_5'] \\
@VVV                  @Vc'VV \\
k[u_4',u_5']  @>d'>>  k[u_5']
\end{CD}$$
Here $\mathcal{O}_Y$ is $k[u_3',u_4',u_5']/(u_3'u_4')$.  On the other charts, no gluing takes place and the variety remains smooth.  Thus $Y$ has only double normal crossing singularities and is semismooth.  

Globally, $Y$ is the pushout 
$$\begin{CD}
\tilde{\mathbb{A}}^1 @>>> Y_2 \\
@VVV                      @VVV \\
Y_1                  @>>> Y
\end{CD}$$
Here $Y_1$ and $Y_2$ are the resolutions of the components $X_1 := (u_0=u_3=0)$ and $X_2 := (u_2 = u_4 =0)$ of $X$, and $\tilde{\mathbb{A}}^1$ is the birational transform of the intersection of $X_1$ and $X_2$.  The universal property of the pushout implies that $f: Y \rightarrow X$ exists.  The open sets of $X$ can be identified with pairs of open sets that have the same pullback to $\mathbb{A}^1$.  Then $Y \rightarrow X$ is an isomorphism over $(U_1,U_2)$, where $U_i$ is the complement of the origin in $X_i$.  This open set clearly has finite complement, so (iv) of the definition is satisfied.  By construction, the singular locus of $Y$ is the birational transform $\mathbb{A}^1$, which is not exceptional.

Finally, we need to verify that $f$ is proper.  We prove this in the following lemma.\qed

\begin{lemma}  Let $X$ be the universal pushout of closed subvarieties $X_1, X_2$ glued along a common subvariety $Z$:
$$\begin{CD}
Z    @>>>  X_2 \\
@VVV       @Vi_2VV \\
X_1  @>i_1>> X
\end{CD}$$
Suppose we are given proper morphisms $g_1: Y_1 \rightarrow X_1$ and $g_2: Y_2 \rightarrow X_2$, and that $Y$ is obtained by gluing $Y_1$ and $Y_2$ along a common subvariety $W$, where $W \rightarrow Z$ is also proper.  Then the induced morphism $f: Y \rightarrow X$ is proper.
\proof  We use the valuative criterion of properness in (\cite{Hart87}, II.4.7).  Suppose that we have a commutative diagram
$$\begin{CD}
\textnormal{Spec } K  @>a>>   Y \\
@VjVV                         @VfVV \\
\textnormal{Spec } R  @>b>>   X 
\end{CD}$$
where $R$ is a valuation ring with quotient field $K$.  Then the morphism $a$ is determined by giving a point $y \in Y$ and an inclusion of residue fields $k(y) \hookrightarrow K$.  Since $Y$ is the union of closed subvarieties $Y_1$ and $Y_2$, and since the residue field of a point of $Y_i$ is the same regardless of whether it is computed in $Y_i$ or $Y$, we see that $a$ lifts to one of the components of $Y$, say $Y_1$.  The composition $Y_1 \hookrightarrow Y \stackrel{f}{\rightarrow} X$ equals $Y_1 \stackrel{g_1}{\rightarrow} X_1 \stackrel{i_1}{\rightarrow} X$, by construction.  Since $g_1$ is proper and so is $i_1$, so is their composition.  Then the valuative criterion of properness implies that there exists a (unique) morphism Spec $R \rightarrow Y_1$ that commutes with the other morphisms.  Composing with $Y_1 \hookrightarrow Y$ gives a morphism that completes the commutative diagram.

We need this morphism to be unique.  However, this follows from the uniqueness of the morphism into $Y_1$.  For, a morphism from Spec $R$ into $Y$ is determined by a pair of points $y_0$ and $y_1$, where $y_1 \in \{ y_0 \}^-$, such that $R$ dominates the local ring of $y_1$ on the subscheme $\{ y_0 \}^-$.  This local ring is the same even if we view $y_0$ as a point of $Y$.  Since we already know that the generic point of Spec $R$ maps into $Y_1$, so must $R$ because $Y_1$ is closed.  So any morphism of Spec $R$ into $Y$ is completely determined by where it sends its generic point, and we conclude that $f$ is proper.\qed
\end{lemma}

\begin{example}  The quotient of the pinch point $(x^2-y^2z = 0)$ by the cyclic group action $x \mapsto x, y \mapsto \epsilon^2y, z \mapsto \epsilon^2z$, where $\epsilon$ generates the cyclic group $\mathbb{Z}/3$.
\end{example}

We construct a semiresolution as before.  The coordinate ring of the quotient is the residue class ring of the ring of invariants $k[x,y^3,y^2z,yz^2,z^3]$ by the intersection with the ideal $(x^2-y^2z)$.  The residue class ring is therefore $k[st,s^3,st^4,t^6]$.  

We can obtain the pinch point by gluing along a quotient map.  To obtain the coordinate ring above, we may also take the ring of invariants of the normalization of the pinch point, and then glue under the quotient map.  Writing the normalization as $k[s,t]$, the group action sends $s \mapsto \epsilon^2s, t \mapsto \epsilon t$.  The ring of invariants is $k[s^3,st,t^3]$.  If we glue by the quotient map, we get the above coordinate ring of $X$.

To find a semiresolution of $X$, we find a resolution of singularities of its normalization, given by $k[s^3,st,t^3]$, keeping track of the double locus defined by $(s^3,st)$.  We then glue along the birational transform of the quotient map to get the semiresolution of $X$.  The normalization has a unique singular point at the origin.  When we blow it up, we obtain three smooth charts, only one of which contains the birational transform of the double locus.  

In particular, writing the normalization as $k[u,v,w]/(v^3-uw)$, the important chart is given by $k[v',w']$.  The birational transform of the double locus is given by $(v'=0)$.  Thus we alter this chart by gluing along the inclusion $k[w'^2] \rightarrow k[w']$.  The fact that we obtain a semiresolution of $X$ this way follows as in the previous example.  Note that the double locus occurs only in this glued chart, as given by the ideal $(v',v'w')$, whereas the exceptional divisor on this chart is $(w'=0)$.  This implies that the semiresolution is strong.
\qed

\noindent
\linebreak
It follows from the definition and Lemma 3.2 that a semiresolution induces an isomorphism (of $\mathcal{O}_X$-modules) of structure sheaves.  We used the $S_2$ property in that proof.  It is worthwhile to note that $f_*\mathcal{O}_Y = \mathcal{O}_X$ also follows from the seminormality of $X$ and $Y$, provided that $f$ has connected fibers.  This may be viewed as a nonnormal version of Zariski's Main Theorem (\cite{Hart87}, III.11.4).  One immediate consequence of that result is that for a normal variety, we get the condition on structure sheaves for any resolution of singularities.  We have the following result in this direction.  Stated for projective morphisms, it is a partial converse to (\cite{Hart87}, III.11.3). 

\begin{lemma}  Suppose $f: Y \rightarrow X$ is a proper and birational morphism with connected fibers such that $Y$ is seminormal.  Then $f_*\mathcal{O}_Y$ is the seminormalization of $\mathcal{O}_X$ in their common sheaf of total quotient rings.
\proof  We consider the Stein factorization of $f$:  $$Y \stackrel{g}{\rightarrow} \textnormal{\textbf{Spec}}_X f_*\mathcal{O}_Y \stackrel{h}{\rightarrow} X.$$  Here $h$ is finite, hence proper.  Then $g$ is also proper by (\cite{Hart87}, II.4.8).  We also know that $h$ is birational.  To see this, suppose that $f: V \rightarrow U$ is an isomorphism.  Then $V$ need not equal the preimage of $U$.  However, the induced morphism $f^{-1}(U) \rightarrow U$ is proper.  The set $f^{-1}(U) \cap V^c$ is closed in $f^{-1}(U)$, so properness implies that its image is closed in $U$.  If we let $$U' = U - f(f^{-1}(U) \cap V^c),$$ then $U'$ is an open subset of $U$ whose preimage lies in $V$.  Then $f^{-1}(U') \rightarrow U'$ is an isomorphism.  In particular, it is an affine morphism, so $h$ is an isomorphism over $U'$.  Since $f$ itself is also birational, we conclude that $g$ is proper and birational.  Then $g$ is surjective.

Since $f$ has connected fibers, so must $h$.  Therefore, $h$ is injective.  Being finite and birational, it is also surjective.  Let $q \in \textnormal{\textbf{Spec}}_X f_*\mathcal{O}_Y$ be the unique point lying over $p \in X$.  The the field extension $k(p) \hookrightarrow k(q)$ is an isomorphism.  In fact, $h$ induces a finite morphism $\{ q \}^- \rightarrow \{ p \}^-$.  By (\cite{LV81}, 2.1), the degree of the field extension $[k(q):k(p)]$ equals the number of points in a generic fiber.  Note that the characteristic of both residue fields is zero, since we are working over a base field of characteristic zero.  Thus $h$ is locally a subintegral extension.  We are done if we can show that the source is seminormal.

We claim that $f_*\mathcal{O}_Y$ is a sheaf of seminormal rings.  The point is that it does not require that we even consider the normalization.  If $A \hookrightarrow B$ is the inclusion of a ring $A$ in its (finite) normalization, then $A$ is seminormal if and only if $b^2, b^3 \in A$ implies $b \in A$.  Equivalently, for every pair of elements $c, d \in A$ such that $c^2 = d^3$, there exists a unique $a \in A$ such that $a^2 =d, a^3 = c$.  The uniqueness guarantees that we can glue sections of $f_*\mathcal{O}_Y$ together to obtain such an element for any ring of sections $\Gamma(U, f_*\mathcal{O}_Y)$.  This finishes the proof.\qed
\end{lemma}

We need to be able to use an arbitrary semiresolution of a variety $X$ to classify its singularities.  In the normal case, one usually shows independence of the resolution by using a single resolution that dominates two given resolutions.  That way, one is able to assume that one resolution factors through the other, and often the proof follows from the fact that the given property is invariant under composition of morphisms.  Thus we state the following.

\begin{proposition}  Let $X$ be a variety, and suppose we are given two semiresolutions $f_1: Y_1 \rightarrow X$ and $f_2: Y_2 \rightarrow X$.  Then there exists a commutative diagram 
$$\begin{CD}
Z     @>q>>    Y_2 \\
@VpVV          @Vf_2VV \\
Y_1   @>f_1>>  X
\end{CD}$$
such that $Z \rightarrow X$ is a semiresolution.
\proof  Let $W = Y_1 \times_X Y_2$ be the fiber product.  It is not clear that $W$ is semismooth in codimension one, even if the product is over a field.  However, we may choose an open set $U \subset X$ such that codim$(X-U,X) \geq 2$, $U$ is semismooth, and both $f_1$ and $f_2$ are isomorphisms over $U$.  For instance, we can just take the intersection of the isomorphism loci.  Then the pullback of $U$ in $W$ is isomorphic to $U$.  

As in the remarks following (2.6), there is a proper morphism $g: Z \rightarrow W$ such that $Z$ is semismooth, $g$ is an isomorphism over $U$, and every component of the singular locus of $Z$ maps onto singular component of $U$.  We get $p$ and $q$ by composing with the projection morphisms.  Note that $Z \rightarrow X$ is proper because properness is preserved under base extension and composition.  By construction, the double locus of $Z$ maps onto (an isomorphic copy of) $U$, which is an open set with sufficiently small complement.  Thus $Z \rightarrow X$ is a semiresolution.

We note that $q: Z \rightarrow Y_2$ (and likewise $p$) is naturally an isomorphism over the codimension one points of $Y_1$ (respectively, over $Y_2$).  In fact, it is an isomorphism over the codimension one points in $U$ by construction.  Every other codimension one point corresponds to a discrete valuation ring, since the semiresolution $Y_2 \rightarrow X$ has the property that the nonnormal codimension one points are in the isomorphism locus.  Thus we can use the valuative criterion for properness to enlarge the rational map from $Y_2$ to $Z$ to include these codimension one points.\qed  
\end{proposition}

In some cases, we would like to have the double locus of a semiresolution intersect the exceptional divisor transversally.  We could call such a semiresolution ``strong."  

\begin{example}  The semiresolution in (3.5) is a strong semiresolution. 
\proof  We noted that one of the charts of the blowup is smooth.  The first chart of the blowup is a pinch point $(y'^2-z'^2w' = 0)$, and the exceptional divisor is $(x'=0)$.  The normalization has coordinates $x', z', y'/z'$.  The double locus is $(y'/z' = 0)$, so it is smooth and intersects the exceptional divisor transversally.  Likewise for the third chart.\qed
\end{example} 

\begin{example}  A semiresolution that is not strong.\end{example}
\noindent
To get a trivial example of a semiresolution that is not strong, start with the pinch point defined by $(x^2-y^2z=0)$, and blow up a point on the positive $z$-axis.  The only singular chart is another pinch point, shifted away from the origin.  The exceptional divisor does not intersect the double locus transversally.  In fact, the double locus is a line that intersects a pair of intersecting lines at their intersection point.\qed

\end{section}

\begin{section}{Semirational Singularities}

We now come to the two main results of this paper.  The first is a semismooth version of G-R vanishing, Theorem 4.2.  It is used in the proof that the definition of semirational singularities does not depend on the choice of a  semiresolution.  The original (smooth) vanishing theorem states that for a resolution of singularities $f: Y \rightarrow X$, we have $R^if_*\omega_Y = 0$ for $i>0$.  In fact, the same condition holds when $f$ is a semiresolution.

\begin{definition}  A variety $X$ has semirational singularities if there is a semiresolution $f: Y \rightarrow X$ such that $R^if_*\mathcal{O}_Y = 0$ for $i > 0$.
\end{definition}

We need to verify that this definition is independent of the semiresolution chosen.  The proof uses the semismooth version of G-R vanishing, which is interesting in its own right.  We keep it as a separate proof so that we can call on it later.

\begin{theorem}  Suppose $f: Y \rightarrow X$ is a semiresolution (of a variety $X$ as in (2.5), $f$ having the properties of (2.7) ).  Then $R^if_*\omega_Y = 0$ for $i > 0$.
\proof  
We consider the commutative diagram 
$$\begin{CD}
\overline{Y}  @>\overline{f}>>  \overline{X} \\
@VpVV                           @VqVV \\
Y             @>f>>							X
\end{CD}$$
where $p$ and $q$ are the normalization morphisms.  Let $D$ be the smooth double locus in $Y$ and $C$ its smooth preimage in $\overline{Y}$.  Then $C \rightarrow D$ is a double covering, ramified along the pinch locus.  Finally, denote by $\mathcal{C}$ the conductor ideal sheaf in $\mathcal{O}_Y$ and denote by $\overline{\mathcal{C}}$ the conductor ideal sheaf in $\mathcal{O}_{\overline{Y}}$.  These ideal sheaves determine the divisors $D$ and $C$, respectively.  Moreover, there is an isomorphism $p_*\overline{\mathcal{C}} \cong \mathcal{C}.$

Applying $p_*$ to a short exact sequence of $\mathcal{O}_{\overline{Y}}$-modules gives another short exact sequence, since $p$ is finite.   Then we have the following commutative diagram with exact rows:

$$\begin{CD}
0 @>>> \mathcal{C}    @>>> \mathcal{O}_Y @>>>                 \mathcal{O}_D    @>>> 0 \\
@VVV   @V\cong VV               @VVV							                  @VVV                  @VVV \\
0 @>>> p_*\overline{\mathcal{C}} @>>> p_*\mathcal{O}_{\overline{Y}} @>>> p_*\mathcal{O}_C @>>> 0
\end{CD}$$

Here we use $p_*\mathcal{O}_C$ to mean the pushforward via the restricted morphism $p|_C$.  Since $Y$ is semismooth, it has an invertible dualizing sheaf $\omega_Y$.  We apply the contravariant cohomological functor $\mathcal{E}xt^*_Y(-, \omega_Y)$ to this diagram.  Since the dualizing sheaf is torsion-free and by (\cite{Hart87}, III.6.3), we obtain a commutative diagram with exact rows

$$\begin{CD}
0 @>>>                  p_*\omega_{\overline{Y}}     @>>> p_*\omega_{\overline{Y}}(C) @>>> \mathcal{E}xt^1_Y(p_*\mathcal{O}_C, \omega_Y) @>>> 0 \\
@VVV                    @VVV                              @V=VV
@VVV                                             @VVV \\
0 @>>> \omega_Y   @>>>  p_*\omega_{\overline{Y}}(C)  @>>> \mathcal{E}xt^1_Y(\mathcal{O}_D, \omega_Y) @>>> 0 \\
\end{CD}$$
Here we have used the identity relating the dualizing sheaf on $Y$ to $\omega_{\overline{Y}}$.  Namely, since $Y$ is $S_2$ and $\omega_Y$ is invertible, we have $p^*\omega_Y = \omega_{\overline{Y}}(C)$.  (We can also pull back differential forms for the hypersurface singularities to obtain this formula directly.)  The terms on the right are $$\mathcal{E}xt^1_Y(\mathcal{O}_D, \omega_Y) = \omega_D \textnormal{ and } \mathcal{E}xt^1_Y(p_*\mathcal{O}_C, \omega_Y) = p_*\omega_C,$$ as can be calculated using duality for a finite morphism.  See (\cite{Reid94}, 2.6) for more details on this exact diagram.  (An important part is that duality implies that $p^{!}\omega_Y = \omega_{\overline{Y}}.$)

Consider the first commutative square above.  We obtain this square using the fact that all the relevant morphisms are birational, so that the normalization of $X$ in $\overline{Y}$ is just the normalization of $X$ in its sheaf of total quotient rings.  Notice that $\overline{f}$ is a resolution of singularities.  It is birational because the other three morphisms are birational, and it is proper because its composition with the finite (proper) morphism $q$ is proper (both $p$ and $f$ are by definition proper).  Finally, we have noted before that $\overline{Y}$ is smooth; the normalization of a semismooth variety is smooth, as can be observed directly by finding the normalizations of double normal crossings and pinch points.  (See the remark following Example 3.4.)  Thus $\overline{f}$ is a resolution of singularities, and we can apply G-R vanishing.  

The spectral sequence for a composition of functors $R^if_*R^jp_*(\omega_{\overline{Y}}) \implies R^{i+j}(f \circ p)_*(\omega_{\overline{Y}})$ degenerates when one of the morphisms is finite.  Thus, going around the above commutative square in opposite directions gives $R^if_*p_*\omega_{\overline{Y}} = q_*R^i\overline{f}_*\omega_{\overline{Y}}.$  These are 0 for $i > 0$ by G-R vanishing for the resolution of singularities $\overline{f}$.  Likewise, the double locus on both $Y$ and $\overline{Y}$ is a smooth variety.  By construction of the semiresolution, the double locus of $Y$ maps birationally onto its image in $X$.  Then the same is true upstairs on $\overline{Y}$.  Again by G-R vanishing and the above spectral sequence argument, $p_*\omega_C$ and $\omega_D$ both have no nonzero higher direct images.

Now we apply $R^if_*$ to the commutative diagram above.  By exactness and the just mentioned vanishing, it follows that the middle term in both rows $p_*\omega_{\overline{Y}}(C)$ has no nonzero higher direct images.  Thus, in the long exact sequence in relative cohomology associated to the bottom row, all terms $R^if_*\omega_Y$ are zero for $i > 1$.  To show that $R^1f_*\omega_Y = 0$, it suffices to show that $f_*p_*\omega_{\overline{Y}}(C) \rightarrow f_*\omega_D$ is surjective.  Since applying $f_*$ to the top row leaves the short exact sequence exact, it is enough to show that $f_*p_*\omega_C \rightarrow f_*\omega_D$ is surjective.

The morphism $p: C \rightarrow D$ is finite, by construction, and flat because $C$ and $D$ are smooth (see (\cite{Hart87}, III.Ex.9.3)).  Thus $p_*\mathcal{O}_C$ is a locally free $\mathcal{O}_D$-module.  There exists a trace morphism $p_*\mathcal{O}_C \rightarrow \mathcal{O}_D$ obtained locally by taking the trace of an element in the free $\mathcal{O}_D$-module $p_*\mathcal{O}_C$ with respect to a suitable basis.  We may normalize this by dividing by the degree of $p$ (note we are in characteristic zero), and then the resulting composition $\mathcal{O}_D \rightarrow p_*\mathcal{O}_C \rightarrow \mathcal{O}_D$ is an identity of $\mathcal{O}_D$-modules.  Applying first $\mathcal{H}om_D(-, \omega_D)$ and then $f_*$ still gives an identity, and in particular the second morphism $f_*p_*\omega_C \rightarrow f_*\omega_D$ is surjective.\qed
\end{theorem}
\begin{theorem}  The condition of Definition 4.1 is independent of the semiresolution chosen.
\proof  For simplicity, we prove this only for projective varieties.  By (3.10), we can restrict to the case where one semiresolution dominates the other.  In this case, we use the spectral sequence for a composition of morphisms $h = f \circ g$, namely $R^if_*(R^jg_*\mathcal{O}_Z) \Rightarrow R^{i+j}h_*\mathcal{O}_Z$.  If we can show that $R^ig_*\mathcal{O}_Z = 0$ for $i > 0$ whenever $g : Z \rightarrow Y$ is a semiresolution of a semismooth variety $Y$, then the spectral sequence degenerates, and we get the desired independence by comparing any two morphisms with a common composition.

We want to use G-R vanishing and Serre duality to show that $R^ig_*\mathcal{O}_Z$ vanishes for $i > 0$.  Suppose that $Y$ is projective.  Choose a very ample divisor $\mathcal{L}$ and an integer $N$ large enough that $R^ig_*\mathcal{O}_Z \otimes \mathcal{L}^N$ is generated by global sections and has no nonzero higher cohomology.  Then it is enough to show that $H^0(Y,R^ig_*\mathcal{O}_Z \otimes \mathcal{L}^N) = 0$.  There is a spectral sequence $H^i(Y,R^jg_*(g^*\mathcal{L}^N)) \Rightarrow H^{i+j}(Z,g^*\mathcal{L}^N)$.  Using the projection formula and the fact that $R^ig_*\mathcal{O}_Z \otimes \mathcal{L}^N$ has no higher cohomology, we see that the spectral sequence degenerates to $H^0(Y, R^ig_*(g^*\mathcal{L}^N)) = H^i(Z, g^*\mathcal{L}^N)$.  Thus it suffices to show that $H^i(Z,g^*\mathcal{L}^N) = 0$ for $i>0$.  Using Serre duality for a projective Cohen-Macaulay scheme, this is equivalent to $H^i(Z, \omega_Z \otimes g^*\mathcal{L}^{-N}) = 0$ for $i < n$, $n = \textnormal{dim }Y$.  

Equivalently, we must show that $H^i(Y,g_*\omega_Z \otimes \mathcal{L}^{-N}) = 0$ for $i< n$.  This is because the spectral sequence $H^i(Y, R^jg_*(\omega_Z \otimes g^*\mathcal{L}^{-N})) \Rightarrow H^{i+j}(Z, \omega_Z \otimes g^*\mathcal{L}^{-N})$ degenerates, by G-R vanishing.  We claim that $g_*\omega_Z = \omega_Y$.  Then the required vanishing is $H^i(Y,\omega_Y \otimes \mathcal{L}^{-N}) = 0$ for $i < n$, which by Serre duality is $H^i(Y,\mathcal{L}^N) = 0$ for $i > 0$.  This is just Serre vanishing when $N$ is sufficiently large.  Note that we are free to apply Serre duality, since a semismooth variety is Cohen-Macaulay.  In fact, its only singularities are hypersurface singularities, so it is Gorenstein.

We need to prove that $g_*\omega_Z = \omega_Y$ for a semiresolution $g$ of semismooth varieties.  First, we observe that $g_*\omega_Z \hookrightarrow (g_*\omega_Z)^{\vee\vee}$ is an inclusion since the pushforward is torsion-free.  The reflexive sheaves $(g_*\omega_Z)^{\vee\vee}$ and $\omega_Y$ are naturally isomorphic in codimension one, since the image of the exceptional divisor has codimension at least two.  Thus there is a natural injection $g_*\omega_Z \hookrightarrow \omega_Y$.  

Next, we look at the expression $\omega_Z = g^*\omega_Y \otimes \mathcal{O}_Z(E)$.  Since $Y$ is semismooth, it has semi canonical singularities.  Thus, as in the smooth case, the exceptional divisor appears with nonnegative coefficients; in other words, there is a morphism $g^*\omega_Y \rightarrow g^*\omega_Y(E) = \omega_Z$.  If we push this forward, we obtain a morphism $\omega_Y \hookrightarrow g_*\omega_Z$, where injectivity follows from reflexivity as above.  To conclude, we have inclusions $\omega_Y \hookrightarrow g_*\omega_Z \hookrightarrow \omega_Y$ whose composition is the identity in codimension one, thus everywhere.  (Reflexive sheaves have property $S_2$.  This is one of the properties of varieties with properties $G_1$ and $S_2$; namely, taking the reflexive hull is an $S_2$-ification.  Considering the cokernel of this injection, any associated point has codimension at least two.  The exact sequence in local cohomology then shows that there can be no associated points, so the cokernel is zero.)  In particular, $g_*\omega_Z \hookrightarrow \omega_Y$ is surjective as well.
\qed
\end{theorem}

Supposing for a moment that we did not require seminormality of $X$ as a base condition, it follows from the definition that $X$ has semirational singularities only if $X$ is seminormal.  In fact, the rings $\Gamma(V, f_*\mathcal{O}_Y)$ are seminormal (see the end of the proof of (3.9)), hence so are the rings $\Gamma(V, \mathcal{O}_X)$ for every open set $V \subset X$.  Seminormality being a local condition, $X$ is therefore seminormal.  Of course, we are using the fact that $f_*\mathcal{O}_Y = \mathcal{O}_X$, which is a consequence of (3.2); as stated in (2.8), we can always choose our semiresolution to be an isomorphism over codimension one points of $X$.

\begin{corollary}  A variety with semirational singularities is Cohen-Macaulay.  
\proof We can just repeat the argument used above.  Suppose that $X$ is projective.  As in (\cite{Hart87}, III.7.6), it suffices to show that $H^i(X, \mathcal{L}^{-N}) = 0$, $i < n$, for an ample sheaf $\mathcal{L}$ and $N$ large.  Semirationality implies that $H^i(X, \mathcal{L}^{-N}) \cong H^i(Y, f^*\mathcal{L}^{-N})$, and then Serre duality on the Cohen-Macaulay variety $Y$ shows that this is isomorphic to $H^{n-i}(Y, \omega_Y \otimes f^*\mathcal{L}^N)'$.  Using semismooth G-R vanishing, this cohomology is isomorphic to $H^{n-i}(X, f_*\omega_Y \otimes \mathcal{L}^N)$.  Finally, Serre vanishing for the coherent sheaf $f_*\omega_Y$ implies that this last cohomology is zero for $n-i > 0$ and $N$ sufficiently large.\qed
\end{corollary}

\begin{corollary}  Semirational singularities are Du Bois.
\proof  This follows from the results in \cite{KSS09} and (\cite{Kov99}, 2.4). In other words, semismooth varieties have Du Bois singularities, since they are Cohen-Macaulay and trivially semi log canonical.  Then the quasi-isomorphism $\mathcal{O}_X \cong_{qis} Rf_*\mathcal{O}_Y$ implies that $X$ also has Du Bois singularities.\qed
\end{corollary}

We have a large class of nonnormal varieties without semirational singularities given by the varieties that are not seminormal.  Even in the restricted case where $X$ is seminormal with $G_1$ and $S_2$, there are plenty of varieties that do not have semirational singularities.  In fact, there are at least as many of these as there are normal varieties that are not Cohen-Macaulay, as the following proof indicates.

\begin{example}  Consider the cone over a normally embedded abelian surface.  Then by gluing with another variety along an appropriate open set, we get a nonnormal variety that is $G_1, S_2$, and seminormal and does not have semirational singularities.
\proof  By (\cite{GT80}, 3.5), any integral normal projective variety $Y$ is birationally equivalent to a seminormal hypersurface $X \subset \mathbb{P}_k^{r+1}$, where dim $Y = r$.  A hypersurface is Gorenstein, hence trivially $G_1$ and $S_2$.  If we glue $X$ to $Y$ along the open set giving a birational equivalence, then all the local properties of the two varieties are preserved.  So the resulting variety is $G_1$, $S_2$, and seminormal, since $X$ and $Y$ both have these properties.  

The cone over a normally-embedded abelian surface (projectively embedded so that the homogeneous coordinate ring is integrally closed) is an example of a normal projective variety that is $S_2$ but not $S_3$ (note that its dimension is 3).  In particular, the cone is not Cohen-Macaulay at the vertex.  So the above construction gives a variety with $G_1$, $S_2$, and seminormality, and which does not have semirational singularities, because these are always Cohen-Macaulay.\qed
\end{example}

\end{section}
\begin{section}{A Note on Non-$\mathbb{Q}$-Gorenstein Varieties}
Recall that a variety is Gorenstein if it is Cohen-Macaulay and has an invertible dualizing sheaf.  It was first intended that a variety should be called $\mathbb{Q}$-Gorenstein if some tensor power of the dualizing sheaf is invertible and the variety is Cohen-Macaulay.  Lately, however, the notion of $\mathbb{Q}$-Gorenstein does not include Cohen-Macaulayness, and so the appellation has become somewhat misleading.  In any case, it should be stated whether one is dealing with Cohen-Macaulay varieties when discussing the  $\mathbb{Q}$-Gorenstein condition.  Note that we almost always assume at least $G_1$, so the dualizing sheaf is actually invertible in codimension one, and that this is even stronger than $S_1$, which is true for reduced varieties.

To classify the singularities of $X$, one uses a resolution (or a semiresolution, when $X$ is nonnormal) $f: Y \rightarrow X$, and relates the canonical divisors by a formula $$K_Y = f^*K_X + \Sigma a_iE_i,$$ where the $E_i$ are the irreducible components of the exceptional divisor.  Depending on whether the $a_i$ are at least -1, or greater than -1, or at least 0, or greater than 0, one obtains definitions of all the pertinent singularities in the birational classification of $X$.

This formula assumes that we can make sense of $K_X$, hence there should be a $G_1$ and $S_2$ hypothesis.  These are the conditions under which a generalized theory of divisors can be worked out.  We also need to be able to pull back $K_X$, and it is here that we use the $\mathbb{Q}$-Gorenstein hypothesis.  The pullback of $K_X$ is a $\mathbb{Q}$-divisor $\frac{1}{m}f^*({mK_X}),$ where $m$ is such that $mK_X$ is Cartier.

In general, it is possible to pull back Cartier divisors (those divisors that correspond to invertible sheaves).  Pulling back arbitrary Weil divisors is somewhat delicate; for instance, the group laws might not be respected, or the intersection numbers might not be well-behaved.  

In their paper ``Singularities on Normal Varieties," de Fernex and Hacon describe a method of overcoming the need for a $\mathbb{Q}$-Gorenstein hypothesis, at least for normal varieties.  See \cite{dFH09}.  It turns out to be the case that these same results hold for nonnormal varieties, using semiresolutions for nonnormal varieties in place of resolutions for normal varieties.  Thus, it is possible to fill in the missing piece of the following diagram, illustrating those classes of singularities that are presently understood:  
$$\begin{array}{ccc}
\textnormal{Normal, $\mathbb{Q}$-Gorenstein} & - & \textnormal{Normal, non-$\mathbb{Q}$-Gorenstein} \\
\vert &  & \vert \\
\textnormal{Nonnormal, $\mathbb{Q}$-Gorenstein} & - & ? \\
\end{array}$$
This has been done by the author in his dissertation, ``Singularities on Nonnormal Varieties."

\nocite{Art70}
\nocite{Kol13}
\nocite{Hart94}
\nocite{Hart87}
\nocite{BH98}
\nocite{GT80}
\nocite{Kov99}
\nocite{LV81}
\nocite{Reid94}
\nocite{dFH09}
\nocite{KSS09}
\nocite{Sch09}

\end{section}

\bibliographystyle{plain}
\bibliography{paper}

\end{document}